\documentclass[10pt]{amsart}
\usepackage{graphicx, color}
\usepackage{amssymb}
\usepackage{amsmath}
\usepackage[T1]{fontenc}
\usepackage[displaymath, mathlines, running]{lineno}
\newtheorem{theorem}{Theorem}[section]
\newtheorem{lemma}[theorem]{Lemma}

\newtheorem{proposition}[theorem]{Proposition}

\newenvironment{proof of main prop}{{\bf Proof of Proposition \ref{main prop}.}}{\hfill\fbox{}\par\vspace{.2cm}}
\newenvironment{proof of main thm}{{\bf Proof of Theorem \ref{main thm}.}}{\hfill\fbox{}\par\vspace{.2cm}}
\newenvironment{proof of prop 3.2}{{\bf Proof of Proposition \ref{prop 3.2}.}}{\hfill\fbox{}\par\vspace{.2cm}}
\newenvironment{proof of prop 3.3}{{\bf Proof of Proposition \ref{prop 3.3}.}}{\hfill\fbox{}\par\vspace{.2cm}}
\newenvironment{proof of prop 5.1}{{\bf Proof of Proposition \ref{prop 5.1}.}}{\hfill\fbox{}\par\vspace{.2cm}}
\numberwithin{equation}{section}

\def\charf {\mbox{{\text 1}\kern-.24em {\text l}}}

\def\bea{\begin{eqnarray*}}
\def\eea{\end{eqnarray*}}
\def\be{\begin{eqnarray}}
\def\ee{\end{eqnarray}}

\begin{document}

\title[Modified Hawking mass and rigidity]{Modified Hawking mass and rigidity of three-manifolds with boundary}


\author{Jihyeon Lee}
\address{Center for Geometry and Physics, Institute for Basic Science (POSTECH Campus), 79 Jigok-ro 127beon-gil, Nam-gu, Pohang, Gyeongbuk, Korea 37673}
\email{jihyeonnie@ibs.re.kr}
\author{Sanghun Lee}
\address{Department of Mathematics and Institute of Mathematical Science, Pusan National University, Busan 46241, Korea}
\email{kazauye@pusan.ac.kr}




\date{\today}

\dedicatory{}

\begin{abstract}
In this paper, we prove a rigidity result for three-dimensional Riemannian manifolds with boundary, under the assumption that a free boundary minimal two-disk, which locally maximizes a modified Hawking mass, is embedded in a $3$-dimensional Riemannian manifold with negative scalar curvature and mean convex boundary. First, we establish area estimates for free boundary strictly stable two-disks. Finally, we show that the $3$-dimensional Riemannian manifold with boundary is locally isometric to the half anti-de Sitter–Schwarzschild manifold.
\end{abstract}

\maketitle
\section{Introduction}

In the study of differential geometry, one of the interesting and fundamental themes lies in understanding the relation between the curvature and topology of a Riemannian manifold. In the two-dimensional case, Schoen and Yau \cite{SC3} demonstrated that any closed surface $\Sigma$ which minimizes area in a three-dimensional Riemannian manifold with positive scalar curvature must be homeomorphic to either $\mathbb{S}^{2}$ or $\mathbb{RP}^{2}$. Motivated by this result, Cai and Galloway \cite{CG} showed that if a three-dimensional Riemannian manifold has nonnegative scalar curvature and contains a closed embedded two-torus $\Sigma$ that is locally area-minimizing, then $\Sigma$ must be flat, and the ambient geometry near $\Sigma$ is locally isometric to the product $\Sigma \times (-\epsilon, \epsilon)$. This rigidity phenomenon was further explored by Bray, Brendle and Neves \cite{BBN}, who extended the result to settings with strictly positive scalar curvature, and by Nunes \cite{NU}, who addressed the case of negative scalar curvature. An alternative approach to this line of work was later provided by Micallef and Moraru \cite{MM}.

M\'{a}ximo and Nunes \cite{MN} recently established an analogous rigidity result in the setting of the de Sitter-Schwarzschild manifold by replacing the area-minimizing assumption with the condition that the surface is a strictly stable minimal two-sphere which maximizes the Hawking mass. Here, we recall the \textit{Hawking mass}, a significant candidate for quasi-local mass in general relativity. For a closed surface $\Sigma \subset M^{3}$, the Hawking mass is defined by
\begin{align*}
m_{H}(\Sigma) = \sqrt{\frac{A(\Sigma)}{16\pi}}\left(\frac{\chi(\Sigma)}{2} - \frac{1}{16\pi}\int_{\Sigma}\left(H^{2} + \frac{2}{3}\inf_{M}R^{M}\right) dv\right),
\end{align*}
where $A(\Sigma)$ is the area of $\Sigma$, $\chi(\Sigma)$ denotes its Euler characteristic, $H$ is the mean curvature of $\Sigma$, and $R^M$ is the scalar curvature of the ambient manifold $M$. The Hawking mass plays a crucial role in geometric analysis, appearing in the proof of the Penrose inequality for asymptotically flat three-manifolds (see \cite{BR, HI}). In particular, the Hawking mass coincides with the ADM mass precisely when the three-manifold is isometric to one-half of the Schwarzschild manifold of $\mathbb{R}^{3}\backslash\{0\}$.

On the other hand, Barros, Batista and Cruz \cite{BBC} established a comparable rigidity result in the setting of the anti-de Sitter–Schwarzschild manifold.

An additional interesting investigation concerns rigidity phenomena in Riemannian manifolds with boundary. In this direction, Ambrozio \cite{AM} extended the results of \cite{BBN, CG, MM, NU} to the setting of Riemannian manifolds with boundary. More recently, the work of M\'{a}ximo and Nunes \cite{MN} was further extended to the boundary case by Batista, Lima and Silva \cite{BLS}, using a modified version of the Hawking mass.

The \textit{modified Hawking mass} is defined as follows. For a compact free boundary surface $\Sigma$ in a three-dimensional Riemannian manifold $M$ with boundary $\partial M$, the modified Hawking mass $\tilde{m}_{H}(\Sigma)$ is defined as
\begin{align} \label{mass}
\tilde{m}_{H}(\Sigma) = \sqrt{\frac{A(\Sigma)}{8\pi}}\left(\chi(\Sigma) - \frac{1}{8\pi}\int_{\Sigma}\left(H^{2} + \frac{2}{3}\inf_{M}R^{M}\right)dv\right).
\end{align}
This notion was first introduced by Marquardt in \cite{MA} and has played an important role in the proof of the Penrose inequality for asymptotically flat manifolds with non-compact boundaries (see also \cite{KO}). Unlike the original Hawking mass, which is defined on closed surfaces, the modified Hawking mass is defined in cases where the surface has boundaries. That is, in the definition of the Hawking mass $m_H$, the Euler characteristic $\chi(\Sigma)$ is $2-2\alpha$ where $\alpha$ is the genus, whereas in the modified Hawking mass $\tilde{m}_H$, the Euler characteristic is defined as $2-2\alpha-\beta$ where $\beta$ is the number of boundary components.

In this paper, our goal is to generalize the rigidity result of Barros, Batista and Cruz \cite{BBC} to the case of manifolds with boundary. Before presenting our results, we provide the definition of the model space, the \textit{half anti-de Sitter-Schwarzschild manifold} (see Example 3.4 in \cite{AME} and Remark 4.3 in \cite{DL}). Fix $m > 0$, the half anti-de Sitter-Schwarzschild manifold $(M_{hadss},g_{hadss})$ is defined as the product space $M_{hadss} := (r_{+},\infty) \times \mathbb{S}^{2}_{+}$ equipped with the metric
\begin{align*}
g_{hadss} := \frac{1}{1 + r^{2} - \frac{2m}{r}}dr^{2} + r^{2}g_{\mathbb{S}_{+}^{2}},
\end{align*}
where $\mathbb{S}_{+}^{2}$ is the hemisphere of $\mathbb{S}^{2}$, $g_{\mathbb{S}_{+}^{2}}$ is the induced metric on $\mathbb{S}^{2}_{+}$, and $r_{+}$ is the unique positive root of the function $f(r) := 1 + r^{2} - \frac{2m}{r}$.

By a change of variables with $\frac{1}{\sqrt{f(r)}}dr=ds$, the half anti-de Sitter-Schwarzschild metric can be rewritten in the warped product form as follows:
\begin{align*}
g_{hadss} = ds^{2} + u(s)^{2}g_{\mathbb{S}^{2}_{+}},     
\end{align*}
where the positive warping function $u: (0,\infty) \rightarrow (r_{+},\infty)$ satisfies $u(s) = r$, $u'(s) = \sqrt{1 + u^{2} - \frac{2m}{u}}$, and $u''(s) = u(s) + \frac{m}{u(s)^{2}}$. 

By reflecting the metric $g_{hadss}$, we construct a complete metric on $\mathbb{R} \times \mathbb{S}^{2}_{+}$ with constant scalar curvature $-6$ and a totally geodesic boundary $\mathbb{R} \times \mathbb{S}^{1}$. The associated warping function $u$ satisfies the following second-order nonlinear differential equation:
\begin{align*}
u'' - \frac{3}{2}u - \frac{1}{2}\left(\frac{1 - (u')^{2}}{u}\right) = 0.
\end{align*}
Based on this, we define a one-parameter family of complete metrics $(g_{hadss})_{a} = ds^{2} + u_{a}(s)^{2}g_{\mathbb{S}_{+}^{2}}$, where each metric has constant scalar curvature equal to $-6$ and boundary mean curvature equal to zero. Here, $u_{a}$ denotes a smooth positive function satisfying the initial conditions $u_{a}(0) = a$ and $u_{a}'(0) = 0$ for $a > 1$. Each metric in this family recovers the half anti-de Sitter–Schwarzschild metric on $\mathbb{R} \times \mathbb{S}^{2}_{+}$ as described above.

In our result, we establish a local rigidity theorem for the half anti-de Sitter–Schwarzschild manifold.

\begin{theorem} \label{main thm}
Let $M$ be a three-dimensional Riemannian manifold with boundary $\partial M$ satisfying $\inf_{M}R^{M} = -6$ and $\inf_{\partial M}H^{\partial M} = 0$. Suppose that $\Sigma$ is a properly embedded, two-sided, free boundary strictly stable minimal two-disk that locally maximizes the modified Hawking mass. Then the Gaussian curvature of $\Sigma$ is constant and equal to $\frac{1}{a^{2}}$, and the geodesic curvature of $\partial\Sigma$ vanishes in a neighborhood of $\Sigma$. Moreover, $M$ is locally isometric to the half anti-de Sitter-Schwarzschild metric $\left(\Sigma \times (-\epsilon,\epsilon), (g_{hadss})_{a}\right)$ for some $\epsilon > 0$.
\end{theorem}

Here, the term \textit{free boundary} means that $\Sigma$ meets $\partial M$ orthogonally along $\partial\Sigma$, and the term \textit{proper} indicates that the boundary $\partial\Sigma$ lies entirely within $\partial M$. 

\section{Preliminaries}

Let $M$ be an $n$-dimensional Riemannian manifold with non-empty boundary $\partial M$, and let $\phi: \Sigma \rightarrow M$ be a smooth two-sided immersion of $\Sigma$ into $M$, where $\Sigma$ is a compact, smooth $(n-1)$-dimensional manifold with non-empty boundary $\partial\Sigma$. We say that an immersion $\phi$ is \textit{proper} if $\phi(\Sigma) \cap \partial M = \phi(\partial \Sigma)$. Let $N$ denote a unit normal vector field along $\Sigma$, and $\nu$ denote the outward unit conormal along $\partial \Sigma$ within $\Sigma$. For simplicity, we identify $\Sigma=\phi(\Sigma)$. The second fundamental form of $\Sigma$ in $M$ is defined by $h^{\Sigma}(X,Y) = \langle \nabla_{X}N, Y\rangle$ for any two tangent vector fields $X$ and $Y$ on $\Sigma$, and the mean curvature $H$ is given by the trace of the second fundamental form $h^{\Sigma}$, that is, $H = {\rm tr}(h^{\Sigma})$. 

Let $\Phi: \Sigma \times (-\epsilon, \epsilon) \rightarrow M$, with $\epsilon>0$, be a smooth map. The one-parameter family of maps $\Phi_s: \Sigma \rightarrow M$, defined by $\Phi_s(x) := \Phi(x, s)$, is called a \textit{variation} of $\phi$ if each $\Phi_s$ is a proper immersion for all $s \in (-\epsilon, \epsilon)$ and satisfies $\Phi_0 = \phi$. The \textit{variational vector field} associated with $\Phi$ is the vector field $V: \Sigma \rightarrow TM$ along $\phi(\Sigma)$, defined by $V(x) = \frac{\partial \Phi}{\partial s}(x, 0)$ for $x \in \Sigma$. For a variation $\Sigma_{s}$ of $\Sigma$, the variation vector field $V$ can be written as the normal variation $V = \varphi N$ for some smooth function $\varphi \in C^{\infty}(\Sigma)$. 
 
The area of $\Sigma$ is given by
\begin{align*}
A(\Sigma) = \int_{\Sigma} \, dv.
\end{align*}
The first variation formula for the area of $\Sigma$ under the variation $\Phi$ with a variational vector field $V = \varphi N$ is
\begin{align*}
\frac{d}{ds}\Bigr|_{s = 0}\, A(\Sigma_{s}) = \int_{\Sigma} \varphi H \, dv + \int_{\Sigma}\varphi\langle N, \nu \rangle \, ds,
\end{align*}
where $\nu$ is the unit conormal to $\partial\Sigma$ in $\Sigma$. The hypersurface $\Sigma$ is said to be \textit{stationary} if it is a critical point for the area functional. This means that a stationary hypersurface $\Sigma$ should be minimal, that is, $H = 0$ on $\Sigma$ as well as $\Sigma$ meets $\partial M$ orthogonally along $\partial\Sigma$. Such a hypersurface $\Sigma$ is called a \textit{free boundary minimal hypersurface}. The volume functional $\text{Vol}: (-\epsilon,\epsilon) \rightarrow \mathbb{R}$ is defined by
\begin{align*}
\text{Vol}(s) := \int_{\Sigma \times [0,s]} \Phi^{*}(dV_{\bar{g}}),
\end{align*}
where $dV_{\bar{g}}$ denote the volume element of $M$. We say that the variation $\Phi$ is \textit{volume-preserving} if $\frac{d}{ds}\Bigr|_{s = 0}\text{Vol}(s) = 0$ for all $s \in (-\epsilon,\epsilon)$. Under the volume-preserving variation, the hypersurface $\Sigma$ is stationary if and only if the mean curvature $H$ is constant on $\Sigma$ and $\Sigma$ meets $\partial M$ orthogonally along $\partial\Sigma$. In this case, we refer to $\Sigma$ as a \textit{constant mean curvature (CMC)} free boundary hypersurface. The second variation formula for the area functional $A(\Sigma_s)$ is given by
\begin{align*}
\frac{d^{2}}{ds^{2}}\Bigr|_{s = 0}\, A(\Sigma_{s}) =& -\int_{\Sigma}\varphi\left(\Delta_{\Sigma}\varphi + \left(Ric(N,N) + \vert h^{\Sigma} \vert^{2}\right)\varphi\right) \, dv \\
&+ \int_{\partial\Sigma}\varphi\left(\frac{\partial\varphi}{\partial\nu} - h^{\partial M}(N,N)\varphi\right) \, ds, \\
=& -\int_{\Sigma}\varphi L_{\Sigma}\varphi \, dv + \int_{\partial\Sigma} \varphi\left(\frac{\partial\varphi}{\partial\nu} - h^{\partial M}(N,N)\varphi\right) \, ds,
\end{align*}
where $L_{\Sigma} = \Delta_{\Sigma} + Ric(N,N) + \vert h^{\Sigma} \vert^{2}$ is the \textit{Jacobi operator} of $\Sigma$ and $h^{\partial M}$ is the second fundamental form of $\partial M$ in M.

It is natural to define the associated quadratic form $Q$ of the second variation, which is called the \textit{index form} of $\Sigma$. For smooth functions $\varphi,\psi \in C^{\infty}(\Sigma)$, the index form $Q$ is given by
\begin{align*}
Q(\varphi,\psi) = -\int_{\Sigma}\varphi L_{\Sigma}\psi \, dv + \int_{\partial\Sigma}\varphi\left(\frac{\partial \psi}{\partial\nu} - h^{\partial M}(N,N)\psi\right) \, ds.
\end{align*}
This bilinear form plays a central role in analyzing the stability of $\Sigma$. 

We now consider the corresponding eigenvalue problem:
\begin{align} \label{eigen prop}
\begin{cases}
L_{\Sigma}\varphi + \lambda\varphi = 0 \,\, &{\rm in} \,\, \Sigma  \\
\frac{\partial\varphi}{\partial\nu} = h^{\partial M}(N,N)\varphi \,\, &{\rm on} \,\, \partial\Sigma.
\end{cases}
\end{align}
By classical spectral theory of elliptic partial differential equations, it is clear that the eigenvalues of this problem are given by an increasing sequence: 
\begin{align*}
\lambda_{1} \leq \lambda_{2} \leq \cdots \nearrow \infty.
\end{align*}
We say that the hypersurface $\Sigma$ is \textit{stable} if the associated index form $Q(\varphi,\varphi)$ is nonnegative for all $\varphi \in C^{\infty}(\Sigma)$. If the smallest eigenvalue $\lambda_{1}$ is strictly positive, the hypersurface $\Sigma$ is called \textit{strictly stable}. Moreover, if $\varphi$ is an eigenfunction corresponding to the first eigenvalue $\lambda_{1}$, then the quadratic form satisfies
\begin{align*}
\lambda_{1}(L_{\Sigma})\int_{\Sigma} \varphi^{2} \, dv = Q(\varphi,\varphi)
\end{align*}
(see \cite{BLS, BCM}).

\section{Area estimates and CMC foliations}

In this section, we derive a sharp estimate for the area of the surface $\Sigma$ in terms of the first eigenvalue of its Jacobi operator. This result relies on the assumption that $\Sigma$ is a strictly stable minimal surface that locally maximizes the modified Hawking mass. The estimate serves as a key step toward the construction of a CMC foliation in a neighborhood of $\Sigma$. Before proving the area estimate, we present several preliminary results, including the first and second variation formulae for the modified Hawking mass, which are essential for our analysis.

\begin{proposition} [\cite{BBC, BLS, MN}] (First variation formula for the modified Hawking mass) 
Let $M$ be a three-dimensional Riemannian manifold with boundary $\partial M$ and $R^{M}$ bounded below. Let $\Sigma$ be a properly embedded, two-sided, free boundary surface. For a proper variation, we have
\begin{align*}
\frac{d}{ds}\Bigr|_{s=0}\tilde{m}_{H}(\Sigma_{s}) =& -\frac{2}{(8\pi)^{\frac{3}{2}}A(\Sigma)^{\frac{1}{2}}}\int_{\Sigma}\left(H\Delta_{\Sigma}\varphi\right) \, dv + \frac{A(\Sigma)^{\frac{1}{2}}}{(8\pi)^{\frac{3}{2}}}\int_{\Sigma}\left(\inf_{M}R^{M} - R^{M}\right)H\varphi \, dv \\
&+ \frac{A(\Sigma)^{\frac{1}{2}}}{(8\pi)^{\frac{3}{2}}}\int_{\Sigma}\left(2K^{\Sigma} - \frac{4\pi\chi(\Sigma)}{A(\Sigma)} - \vert h^{\Sigma} \vert^{2} + \frac{1}{2A(\Sigma)}\int_{\Sigma}H^{2} \, dv\right)H\varphi \, dv,
\end{align*}
where $K^{\Sigma}$ is the Gaussian curvature of $\Sigma$.
\end{proposition}

We now present the second variation formula for the modified Hawking mass.

\begin{proposition} [\cite{BBC, BLS, MN}] (Second variation formula for the modified Hawking mass) \label{secondvariation}
Let $M$ be a three-dimensional Riemannian manifold with boundary $\partial M$ and $R^{M}$ bounded below. Let $\Sigma$ be a properly embedded, two-sided, free boundary surface. If $\Sigma$ is a critical point of the modified Hawking mass, then we have
\begin{align*}
\frac{d^{2}}{ds^{2}}\Bigr|_{s=0}\tilde{m}_{H}(\Sigma_{s}) =& -\frac{3\tilde{m}_{H}(\Sigma)}{4A(\Sigma)^{2}}\left(\int_{\Sigma}H^{2}\varphi \, dv \right)^{2} - \frac{2A(\Sigma)^{\frac{1}{2}}}{(8\pi)^{\frac{3}{2}}}\int_{\Sigma}\left((L_{\Sigma}\varphi)^{2} + HL'_{\Sigma}(0)\varphi\right) \, dv \\
&+ \frac{A(\Sigma)^{\frac{1}{2}}}{(8\pi)^{\frac{3}{2}}}\int_{\Sigma}\left(H^{2} + \frac{2}{3}\inf_{M}R^{M}\right)\left(\varphi L_{\Sigma}\varphi - H^{2}\varphi^{2}\right) \, dv \\
&+ \frac{4A(\Sigma)^{\frac{1}{2}}}{(8\pi)^{\frac{3}{2}}}\int_{\Sigma}\left(H^{2}\varphi L_{\Sigma}\varphi\right) \, dv \\
&- \frac{\tilde{m}_{H}(\Sigma)}{2A(\Sigma)}\int_{\Sigma}\left(\varphi L_{\Sigma}\varphi - H^{2}\varphi^{2} + div_{\Sigma}(\nabla_{X}X)\right) \, dv \\
&- \frac{A(\Sigma)^{\frac{1}{2}}}{(8\pi)^{\frac{3}{2}}}\int_{\partial\Sigma}\left(H^{2} + \frac{2}{3}\inf_{M}R^{M}\right)\left(\frac{\partial\varphi}{\partial\nu} - h^{\partial M}(N,N) \varphi \right)\varphi \, ds \\
&+ \frac{\tilde{m}_{H}(\Sigma)}{2A(\Sigma)}\int_{\partial\Sigma}\left(\frac{\partial\varphi}{\partial\nu} - h^{\partial M}(N,N) \varphi \right)\varphi \, ds,
\end{align*}
where $L'_{\Sigma}(0)$ is the first variation of the Jacobi operator $L_{\Sigma}$ in \cite{BLS, MN}.
\end{proposition}
Further details on the first variation of the Jacobi operator $L_{\Sigma}$ can be found in \cite[Proposition 2]{BLS} and \cite[Proposition A.2]{MN}. 

We now present a key area estimate for minimal surfaces under a modified Hawking mass maximization condition. The following proposition provides a precise relation between the area of a free boundary minimal disk and the first eigenvalue of its Jacobi operator, under appropriate curvature assumptions on the ambient manifold.

\begin{proposition} \label{main prop}
Let $M$ be a three-dimensional Riemannian manifold with boundary $\partial M$ satisfying $\inf_{M}R^{M} = - 6$ and $\inf_{\partial M}H^{\partial M} = 0$. Suppose $\Sigma$ is a properly embedded, two-sided, free boundary strictly stable minimal two-disk that locally maximizes the modified Hawking mass. Then the area of $\Sigma$ satisfies the estimate
\begin{align*}
A(\Sigma) = \frac{2\pi}{\lambda_{1}(L_{\Sigma}) - 3},
\end{align*}
where $\lambda_{1}(L_{\Sigma})$ is the first eigenvalue of the Jacobi operator $L_{\Sigma}$. Moreover, we have that 
\begin{itemize}
\item[(i)] $\Sigma$ is totally geodesic,
\item[(ii)] $Ric(N,N) = -\lambda_{1}(L_{\Sigma})$, $R^{M} = -6$, and $K^{\Sigma} = \frac{2\pi}{A(\Sigma)}$ in $\Sigma$,
\item[(iii)] $H^{\partial M} = 0$ and $k_{g} = 0$ along $\partial\Sigma$,
\end{itemize}
where $K^{\Sigma}$ is the Gaussian curvature of $\Sigma$ and $k_{g}$ is the geodesic curvature of $\partial\Sigma$.
\end{proposition}

\begin{proof of main prop}
Since $\Sigma$ is a strictly stable minimal surface, we have
\begin{align*}
Q(\varphi,\varphi) \geq \lambda_{1}(L_{\Sigma})\int_{\Sigma}\varphi^{2} \, dv
\end{align*}
for some function $\varphi \in C^{\infty}(\Sigma)$. Taking $\varphi \equiv 1$ on $\Sigma$ yields that
\begin{align*}
\lambda_{1}(L_{\Sigma})A(\Sigma) 
& \leq -\int_{\Sigma}\left(Ric(N,N) + \vert h^{\Sigma} \vert^{2}\right) \, dv - \int_{\partial\Sigma} h^{\partial M}(N,N) \, ds \\
& \leq \int_{\Sigma} \left(-\frac{R^{M}}{2} + K^{\Sigma} - \frac{\vert h^{\Sigma} \vert^{2}}{2}\right) dv - \int_{\partial\Sigma} \left(H^{\partial M} - k_{g}\right) ds \\
&\leq 3A(\Sigma) + 2\pi\chi(\Sigma) \\
&= 3A(\Sigma) + 2\pi
\end{align*}
where the second inequality holds by the Gauss equation and $h^{\partial M}(N,N) = H^{\partial M} - k_{g}$, and the third inequality holds by the assumption for $R^{M}$ and $H^{\partial M}$ and the Gauss-Bonnet theorem. 
Therefore we obtain
\begin{align} \label{area esi}
A(\Sigma) \leq \frac{2\pi}{\lambda_{1}(L_{\Sigma}) - 3}.
\end{align}
If the equality holds, then all the above inequalities become equalities. This implies that $Q(1,1) = \lambda_{1}(L_{\Sigma})A(\Sigma)$. Therefore $\Sigma$ is totally geodesic, $R^{M} = -6$, and $H^{\partial M} = 0$ in $\Sigma$. In order to obtain more properties, we define a new functional $F$ as follows:
\begin{align*}
F(\varphi,\psi) := Q(\varphi,\psi) - \lambda_{1}(L_{\Sigma})\int_{\Sigma}\varphi\psi \, dv,
\end{align*}
for any $\varphi,\psi \in C^{\infty}(\Sigma)$. 
Then $F(\varphi,\varphi) \geq 0$ and $F(1,1) = 0$ for any $\varphi \in C^{\infty}(\Sigma)$. This implies $F(1,\varphi) = 0$ for every $\varphi \in C^{\infty}(\Sigma)$. So we have
\begin{align*}
0 = F(1,\varphi) = - \int_{\Sigma}\left(Ric(N,N) + \lambda_{1}(L_{\Sigma})\right)\varphi \, dv - \int_{\partial\Sigma}h^{\partial M}(N,N)\varphi \, ds.
\end{align*}
Hence, we deduce that $Ric(N,N) = - \lambda_{1}(L_{\Sigma})$, $K^{\Sigma} = \frac{2\pi}{A(\Sigma)}$ in $\Sigma$, and $k_{g} = 0$ along $\partial\Sigma$.

Now we prove the reverse direction of the inequality \eqref{area esi}. Since $\Sigma$ locally maximizes the modified Hawking mass, from Proposition \ref{secondvariation}, we have
\begin{align*}
0 &\geq -\frac{2A(\Sigma)^{\frac{1}{2}}}{(8\pi)^{\frac{3}{2}}}\int_{\Sigma}(L_{\Sigma}\varphi)^{2} \, dv - \left(\frac{4A(\Sigma)^{\frac{1}{2}}}{(8\pi)^{\frac{3}{2}}} + \frac{\tilde{m}_{H}(\Sigma)}{2A(\Sigma)}\right)\int_{\Sigma} \varphi L_{\Sigma}\varphi \, dv \\
&\quad +\left(\frac{4A(\Sigma)^{\frac{1}{2}}}{(8\pi)^{\frac{3}{2}}} + \frac{\tilde{m}_{H}(\Sigma)}{2A(\Sigma)}\right)\int_{\partial\Sigma}\left(\frac{\partial\varphi}{\varphi\nu} - \langle N, \nabla_{N}N^{\partial M} \rangle\varphi\right)\varphi \, ds.
\end{align*}
By the eigenvalue problem \eqref{eigen prop} and the normalized condition $\int_{\Sigma}\varphi^{2} \, dv = 1$, we obtain
\begin{align*}
0 \geq -\frac{2A(\Sigma)^{\frac{1}{2}}}{(8\pi)^{\frac{3}{2}}}\lambda_{1}(L_{\Sigma})^{2} + \frac{4A(\Sigma)^{\frac{1}{2}}}{(8\pi)^{\frac{3}{2}}}\lambda_{1}(L_{\Sigma}) + \frac{\tilde{m}_{H}(\Sigma)}{2A(\Sigma)}\lambda_{1}(L_{\Sigma}).
\end{align*}
Note that the modified Hawking mass $\tilde{m}_H$ is exactly given by $\tilde{m}_{H}(\Sigma) = \frac{A(\Sigma)^{\frac{1}{2}}}{(8\pi)^{\frac{1}{2}}} + \frac{4A(\Sigma)^{\frac{3}{2}}}{(8\pi)^{\frac{3}{2}}}$ when $\inf_{M}R^{M} = -6$ and $\Sigma$ is a minimal two-disk. Combining this with the above inequality, we have
\begin{align*}
\lambda_{1}(L_{\Sigma})^{2} \geq 3\lambda_{1}(L_{\Sigma}) + \frac{2\pi\lambda_{1}(L_{\Sigma})}{A(\Sigma)}.
\end{align*}
Since $\lambda_{1}(L_{\Sigma}) > 0$, we get
\begin{align*}
A(\Sigma) \geq \frac{2\pi}{\lambda_{1}(L_{\Sigma}) - 3}.
\end{align*}
\end{proof of main prop}

The following proposition concerns the construction of CMC foliations near a free boundary strictly stable minimal two-disk in a three-dimensional Riemannian manifold with boundary. The derivation of this result follows similar arguments to those in \cite{BBC, BLS, MN}. For detailed proofs and techniques, we refer the reader to \cite[Proposition 2]{BBC}, \cite[Proposition 5]{BLS} and \cite[Proposition 5.1]{MN}.

\begin{proposition} [\cite{BBC, BLS, MN}] \label{CMCfoliatioin}
 Let $M$ be a three-dimensional Riemannian manifold with boundary $\partial M$. Consider a properly embedded, two-sided, free boundary strictly stable minimal two-disk $\Sigma$. If $\Sigma$ satisfies the following area estimate
\begin{align*}
A(\Sigma) = \frac{2\pi}{\lambda_{1}(L_{\Sigma}) - 3},
\end{align*}
then there exists a positive real number $\epsilon$ and a smooth map $\mu: (-\epsilon,\epsilon) \times \Sigma \rightarrow \mathbb{R}$ satisfying the following two properties:
\begin{itemize}
    \item [(i)] $\mu(0,x) = 0, \quad \frac{\partial\mu}{\partial t}(0,x) = 1, \quad \text{and} \quad \int_{\Sigma}\left(\mu(t,x) - t\right)dv = 0$,
    \item [(ii)] $\Sigma_{t} \, \text{has free boundary CMC two-disk for all} \,\, t \in (-\epsilon,\epsilon)$.
\end{itemize}
\end{proposition}

By applying the above proposition, we can choose positive $\epsilon$ sufficiently small so that the mean curvature function $H(t)$ of $\Sigma_t$ is strictly positive for $ t \in (-\epsilon, 0)$ and strictly negative for $t \in (0, \epsilon)$; that is, $H(t)$ maintains a definite sign on either side of $t = 0$. Let $b_t(x) = \exp_x(\mu(x,t) N_t(x))$ be a smooth variation of the surface $\Sigma$, where $\mu(x,t)$ is a smooth function controlling the normal displacement, \( N_t(x) \) is the unit normal vector field to the deformed surface \( \Sigma_t \), and \( x \in \Sigma \), \( t \in (-\epsilon, \epsilon) \). The \textit{lapse function} \( \rho_t \) is defined by
\begin{align*}
    \rho_t(x) = \left\langle N_t(x), \frac{\partial}{\partial t} b_t(x) \right\rangle,
\end{align*}
which measures the rate at which the deformation moves in the direction of the normal vector at each point. Let $H(t)$ denote the mean curvature of the hypersurface $\Sigma_t$. Since $\{\Sigma_t\}_{t \in (-\epsilon, \epsilon)}$ forms a CMC foliation around the minimal surface $\Sigma = \Sigma_0$ by Proposition \ref{CMCfoliatioin}, the mean curvature function $H(t)$ satisfies the following linearized system:
\begin{align*}
\begin{cases}
H'(t) = L_{\Sigma_t} \rho_t & \text{in } \Sigma_t, \\
\frac{\partial \rho_t}{\partial \nu_t} = h^{\partial M}(N_t, N_t)\rho_t & \text{along } \partial\Sigma_t,
\end{cases}
\end{align*}
where $L_{\Sigma_t}$ is the Jacobi operator associated with $\Sigma_t$, and $\nu_t$ is the outward conormal along $\partial\Sigma_t$. It can be easily checked that at $t = 0$, we have $\rho_0 \equiv 1$, since the variation is normal and unit speed at $\Sigma$, and the derivative of the mean curvature at $t = 0$ satisfies
\begin{align*}
    \left.\frac{d}{dt}\right|_{t=0} H(t) = L_{\Sigma} \rho_0 = L_{\Sigma}(1) = -\lambda_1(L_{\Sigma}) < 0,
\end{align*}
by the definition of the first eigenvalue $\lambda_1$ and the assumption that $\Sigma$ is strictly stable. Therefore, $H(t)$ decreases near $t = 0$, and by continuity, we may choose $\epsilon > 0$ sufficiently small such that $H(t) > 0$ for $t \in (-\epsilon, 0)$, and $H(t) < 0$ for $t \in (0, \epsilon)$, meaning the mean curvature changes sign across \( \Sigma \), consistent with \( \Sigma \) being a local maximum of the modified Hawking mass.

\section{Proof of Theorem \ref{main thm}}

In this section, we present the proof of our main result. In order to support the proof of Theorem \ref{main thm}, we begin by recalling the following lemma.

\begin{lemma} [Lemma 1 in \cite{BLS}] \label{lemma 4.1}
Let $\Sigma_{t}$ be free boundary CMC foliations. Then
\begin{align*}
\int_{\Sigma_{t}}\left(Ric(N_{t},N_{t}) + \vert h^{\Sigma}_{t} \vert^{2}\right)\rho_{t} \, dv &= \bar{\rho}_{t}\int_{\Sigma_{t}}\left(Ric(N_{t},N_{t}) + \vert h^{\Sigma}_{t} \vert^{2}\right) dv + \bar{\rho}_{t}\int_{\Sigma_{t}}\frac{\vert \nabla^{\Sigma_{t}}\rho_{t} \vert^{2}}{\rho^{2}_{t}} \, dv \\
&\quad + H'(t)\theta(t,x) - \int_{\Sigma_{t}}\Delta_{\Sigma_{t}}\rho_{t} \, dv + \bar{\rho}_{t}\int_{\partial\Sigma_{t}}h^{\partial M}(N_{t},N_{t}) \, ds,
\end{align*}
where $\bar{\rho}_{t} = \frac{1}{A(\Sigma_{t})}\int_{\Sigma_{t}}\rho_{t} \, dv$ and $\theta(t,x)$ is a non-positive function.
\end{lemma}

Using Lemma \ref{lemma 4.1} above, we now proceed to prove Theorem \ref{main thm}. 

\begin{proof of main thm}
In $\Sigma_{t}$, the modified Hawking mass is defined as
\begin{align} \label{CMC mass}
\tilde{m}_{H}(\Sigma_{t}) = \frac{A(\Sigma_{t})^{\frac{1}{2}}}{(8\pi)^{\frac{3}{2}}}\left(8\pi - \int_{\Sigma_{t}}\left(H(t)^{2} + \frac{2}{3}\inf_{M}R^{M}\right) dv\right).
\end{align}
If we differentiate the equation (\ref{CMC mass}), we have
\begin{align} \label{first mass}
\frac{d}{dt}\tilde{m}_{H}(\Sigma_{t}) = &-\frac{A(\Sigma_{t})^{-\frac{1}{2}}}{2(8\pi)^{\frac{3}{2}}}\left(8\pi - \int_{\Sigma_{t}}\left(H(t)^{2} + \frac{2}{3}\inf_{M}R^{M}\right)dv\right)\left(\int_{\Sigma_{t}}H(t)\rho_{t} \, dv\right) \\
&- \frac{2A(\Sigma_{t})^{\frac{1}{2}}}{(8\pi)^{\frac{3}{2}}}\int_{\Sigma_{t}}H(t)\left(\Delta_{\Sigma_{t}}\rho_{t} + Ric(N_{t},N_{t})\rho_{t} + \vert h^{\Sigma}_{t} \vert^{2}\rho_{t}\right)dv \nonumber \\
&+ \frac{A(\Sigma_{t})^{\frac{1}{2}}}{(8\pi)^{\frac{3}{2}}}\int_{\Sigma_{t}}H(t)\left(H(t)^{2} + \frac{2}{3}\inf_{M}R^{M}\right)\rho_{t} \, dv \nonumber \\
= &- \frac{A(\Sigma_{t})^{\frac{1}{2}}}{(8\pi)^{\frac{3}{2}}}H(t)\left(4\pi\bar{\rho}_{t} - A(\Sigma_{t})\left(\frac{H(t)^{2}}{2} + \frac{1}{3}\inf_{M}R^{M}\right)\bar{\rho}_{t}\right) \nonumber \\
&- 2\frac{A(\Sigma_{t})^{\frac{1}{2}}}{(8\pi)^{\frac{3}{2}}}H(t)\int_{\Sigma_{t}}\Delta_{\Sigma_{t}}\rho_{t} \, dv + \frac{A(\Sigma_{t})^{\frac{3}{2}}}{(8\pi)^{\frac{3}{2}}}H(t)\left(H(t)^{2} + \frac{2}{3}\inf_{M}R^{M}\right)\bar{\rho}_{t} \nonumber \\
&- 2\frac{A(\Sigma_{t})^{\frac{1}{2}}}{(8\pi)^{\frac{3}{2}}}H(t)\int_{\Sigma_{t}}\rho_{t}\left(Ric(N_{t},N_{t}) + \vert h^{\Sigma}_{t}\vert^{2}\right)dv \nonumber \\
= &-\frac{A(\Sigma_{t})^{\frac{1}{2}}}{(8\pi)^{\frac{3}{2}}}H(t)\left(4\pi\bar{\rho}_{t} - A(\Sigma_{t})\bar{\rho}_{t}\left(\frac{3}{2}H(t)^{2} + \inf_{M}R^{M}\right)\right) \nonumber \\
&-\frac{A(\Sigma_{t})^{\frac{1}{2}}}{(8\pi)^{\frac{3}{2}}}H(t)\left(2\int_{\Sigma_{t}}\Delta_{\Sigma_{t}}\rho_{t} \, dv + 2\int_{\Sigma_{t}}\rho_{t}\left(Ric(N_{t},N_{t}) + \vert h^{\Sigma}_{t} \vert^{2}\right) dv\right). \nonumber
\end{align}
On the other hand, by the Lemma \ref{lemma 4.1} and the Gauss equation, we obtain
\begin{align} \label{eq 4.3}
2\int_{\Sigma_{t}}\rho_{t}\left(Ric(N_{t},N_{t}) + \vert h^{\Sigma}_{t} \vert^{2}\right)dv 
&= \bar{\rho}_{t}\int_{\Sigma_{t}}R^{M} \, dv - 2\bar{\rho}_{t}\int_{\Sigma_{t}}K^{\Sigma_{t}} \, dv + \bar{\rho}_{t}\int_{\Sigma_{t}}H(t)^{2} \, dv  \\
&\quad + \bar{\rho}_{t}\int_{\Sigma_{t}}\vert h^{\Sigma}_{t} \vert^{2} \, dv + 2\bar{\rho}_{t}\int_{\Sigma_{t}}\frac{\vert \nabla^{\Sigma_{t}}\rho_{t} \vert^{2}}{\rho_{t}^{2}} \, dv \nonumber \\
&\quad -2\int_{\Sigma_{t}}\Delta_{\Sigma_{t}}\rho_{t} \, dv + 2H'(t)\theta(t,x) \nonumber \\
&\quad + 2\bar{\rho}_{t}\int_{\partial\Sigma_{t}}H^{\partial M} \, ds - 2\bar{\rho}_{t}\int_{\partial\Sigma_{t}}k_{g,t} \, ds \nonumber \\
&= \bar{\rho}_{t}\int_{\Sigma_{t}}R^{M} \, dv + \bar{\rho}_{t}\int_{\Sigma_{t}}H(t)^{2} \, dv + \bar{\rho}_{t}\int_{\Sigma_{t}}\vert h^{\Sigma}_{t} \vert^{2} \, dv \nonumber \\
&\quad + 2\bar{\rho}_{t}\int_{\Sigma_{t}}\frac{\vert \nabla^{\Sigma_{t}}\rho_{t} \vert^{2}}{\rho_{t}^{2}} \, dv - 2\int_{\Sigma_{t}}\Delta_{\Sigma_{t}}\rho_{t} \, dv \nonumber \\
&\quad + 2H'(t)\theta(t,x) +2\bar{\rho}_{t}\int_{\partial\Sigma_{t}}H^{\partial M} \, ds - 4\pi\bar{\rho}_{t}, \nonumber
\end{align}
where the second equality holds by the Gauss-Bonnet theorem. 
If we substitute the equation \eqref{eq 4.3} into \eqref{first mass}, we get
\begin{align} \label{first mass 2}
\frac{d}{dt}\tilde{m}_{H}(\Sigma_{t}) = &-\frac{A(\Sigma_{t})^{\frac{1}{2}}}{(8\pi)^{\frac{3}{2}}}H(t)\left(\bar{\rho}_{t}\int_{\Sigma_{t}}\left(R^{M} - \inf_{M}R^{M}\right)dv + 2H'(t)\theta(t,x)\right) \\
&- \frac{A(\Sigma_{t})^{\frac{1}{2}}}{(8\pi)^{\frac{3}{2}}}H(t)\left(\bar{\rho}_{t}\int_{\Sigma_{t}}\left(\vert h^{\Sigma}_{t} \vert^{2} - \frac{H(t)^{2}}{2}\right)dv + 2\bar{\rho}_{t}\int_{\Sigma_{t}}\frac{\vert \nabla^{\Sigma_{t}}\rho_{t} \vert^{2}}{\rho_{t}^{2}} \, dv \right) \nonumber \\
&- 2\frac{A(\Sigma_{t})^{\frac{1}{2}}}{(8\pi)^{\frac{3}{2}}}H(t)\bar{\rho}_{t}\int_{\partial\Sigma_{t}}H^{\partial M} \, ds. \nonumber
\end{align}
Since $\rho_{0}(x) = 1$ for all $x \in \Sigma$, we can assume that $\rho_{t}(x) > 0$ for sufficiently small $\epsilon$ and $t \in (-\epsilon,\epsilon)$. We remark that $H(t) > 0$ for $t \in (-\epsilon,0)$ and $H(t) < 0$ for $t \in (0,\epsilon)$.
\begin{itemize}
\item [CASE 1.] $H(t) > 0$ for $t \in (-\epsilon,0)$.
\end{itemize}
With the assumption for $R^{M}$ and $H^{\partial M}$, we obtain
\begin{align*}
-\frac{A(\Sigma_{t})^{\frac{1}{2}}}{(8\pi)^{\frac{3}{2}}}H(t)\bar{\rho}_{t}\int_{\Sigma_{t}}\left(R^{M} - \inf_{M}R^{M}\right)dv \leq 0
\end{align*}
and
\begin{align*}
-2\frac{A(\Sigma_{t})^{\frac{1}{2}}}{(8\pi)^{\frac{3}{2}}}H(t)\bar{\rho}_{t}\int_{\partial\Sigma_{t}}H^{\partial M} \, ds \leq 0.
\end{align*}
Moreover, since $\vert h^{\Sigma}_{t} \vert^{2} \geq \frac{H(t)^{2}}{2}$ and $ \frac{\vert \nabla^{\Sigma_{t}}\rho_{t} \vert^{2}}{\rho_{t}^{2}} \geq 0$, we get
\begin{align*}
- \frac{A(\Sigma_{t})^{\frac{1}{2}}}{(8\pi)^{\frac{3}{2}}}H(t)\left(\bar{\rho}_{t}\int_{\Sigma_{t}}\left(\vert h^{\Sigma}_{t} \vert^{2} - \frac{H(t)^{2}}{2}\right)dv + 2\bar{\rho}_{t}\int_{\Sigma_{t}}\frac{\vert \nabla^{\Sigma_{t}}\rho_{t} \vert^{2}}{\rho_{t}^{2}} \, dv \right) \leq 0.
\end{align*}
Since $H(t) > 0$, $H(0) = 0$, and $\theta(t,x) \leq 0$ for $t \in (-\epsilon, 0)$, we have
\begin{align*}
-2\frac{A(\Sigma_{t})^{\frac{1}{2}}}{(8\pi)^{\frac{3}{2}}}H(t)\left(H'(t)\theta(t,x)\right) \leq 0.
\end{align*}
Therefore, $\frac{d}{dt}\tilde{m}_{H}(\Sigma_{t}) \leq 0$ for all $t \in (-\epsilon, 0)$.
\begin{itemize}
\item [CASE 2.] $H(t) < 0$ for $t \in (0, \epsilon)$.
\end{itemize}
Similar to CASE 1, we can infer $\frac{d}{dt}\tilde{m}_{H}(\Sigma_{t}) \geq 0$ for all $t \in (0,\epsilon)$.

Combining CASE 1 and CASE 2, we deduce
\begin{align*}
\tilde{m}_{H}(\Sigma_{t}) \geq \tilde{m}_{H}(\Sigma)
\end{align*}
for all $t \in (-\epsilon,\epsilon)$. But our assumption is that $\Sigma$ is locally maximizes the modified Hawking mass. So we conclue that $\frac{d}{dt}\tilde{m}_{H}(\Sigma_{t}) = 0$. Hence, this implies that $R^{M} = -6$, $\Sigma_{t}$ is umbilic in $\Sigma_{t}$, and $H^{\partial M} = 0$ along $\partial\Sigma_{t}$. By the equation \eqref{first mass 2}, we have
\begin{align*}
H'(t)\theta(t,x) + \bar{\rho}_{t}\int_{\Sigma_{t}}\frac{\vert\nabla^{\Sigma_{t}} \rho_{t}\vert^{2}}{\rho_{t}^{2}} = 0
\end{align*}
for $t \in (-\epsilon,\epsilon)$. But we know that $H'(t)\theta(t,x) \geq 0$. It follows that $\rho_{t} = 1$ for all $t \in (-\epsilon,\epsilon)$. Therefore, we can choose a small neighborhood of $\Sigma$ such that the metric is given by $g = dt^{2} + g_{\Sigma_{t}}$, where $g_{\Sigma_{t}}$ is the induced metric by the isometry $b(x,t) = {\rm exp}_{x}(tN(x))$. Applying Theorem 3.2 in \cite{HP} gives that
\begin{align*}
\frac{\partial}{\partial t}g_{\Sigma_{t}} &= -2\rho_{t}h^{\Sigma}_{t} \\
&= -H(t)g_{\Sigma_{t}},
\end{align*}
where in the last equality we used $\rho_{t} = 1$, $H(t)$ is constant and $\Sigma_{t}$ is umbilic for all $t \in (-\epsilon,\epsilon)$. We conclude that $g_{\Sigma_{t}} = u_{a}(t)^{2}g_{S^{2}_{+}}$, where $u_{a}(t) = ae^{-\frac{1}{2}\int^{t}_{0}H(s)ds}$ and $a^{2} = \frac{A(\Sigma)}{2\pi}$ for all $t \in (-\epsilon,\epsilon)$. By the uniqueness of the solution of ODE, we have that $M$ is isometric to $(\Sigma \times (-\epsilon,\epsilon),(g_{hadss})_{a})$ in a neighborhood of $\Sigma$.
\end{proof of main thm}

\end{document}